\documentclass{amsart}

\newtheorem{theorem}{Theorem}[section]

\newtheorem{corollary}[theorem]{Corollary}
\newtheorem{os}[theorem]{The Odlyzko-Sch\"onhage algorithm}

\theoremstyle{definition}

\theoremstyle{remark}

\numberwithin{equation}{section}

\begin{document}

\title[A method to compute the zeta function]{An amortized-complexity method to compute the Riemann zeta function}

\author[G.A. Hiary]{Ghaith A. Hiary}
\address{Institute for advanced Study, 1 Einstein Drive, Princeton, NJ, 08540.} 
\email{hiaryg@gmail.com}
\thanks{This material is based upon work supported by the National Science Foundation under agreements No. DMS-0757627 (FRG grant) and No. DMS-0635607.}

\date{}

\begin{abstract}
A practical method to compute the Riemann zeta function is presented. The method can compute $\zeta(1/2+it)$ at any $\lfloor T^{1/4} \rfloor$ points in $[T,T+T^{1/4}]$ using an {\em{average}} time of $T^{1/4+o(1)}$ per point. This is the same complexity as the Odlyzko-Sch\"onhage algorithm over that interval. Although the method far from competes with the Odlyzko-Sch\"onhage algorithm over intervals much longer than $T^{1/4}$, it still has the advantages of being elementary, simple to implement, it does not use the fast Fourier transform or require large amounts of storage space, and its error terms are easy to control. The method has been implemented, and results of timing experiments agree with its theoretical amortized complexity of $T^{1/4+o(1)}$.
\end{abstract}

\maketitle

\section{Introduction}

The Riemann zeta function $\zeta(s)$ can be calculated using the Riemann-Siegel formula. A frequently stated version of that formula on the critical line is: let

\begin{equation} \label{eq:Z(t)}
Z(t) = e^{i\theta(t)} \zeta(1/2+it)\,,\qquad e^{i\theta(t)} = \left(\frac{\Gamma(1/4+it/2)}{\Gamma(1/4-it/2)}\right)^{1/2} \pi^{-it/2}\,.  
\end{equation}

\noindent
The rotation factor $e^{i\theta(t)}$ is chosen so that $Z(t)$ is real. Let $a:=\sqrt{t/(2\pi)}$, $n_1:=\lfloor a\rfloor$ be the integer part of $a$, and $\rho:= \{a\}=a-\lfloor a\rfloor$ be the fractional part of $a$. Then for $t> 2\pi$,  

\begin{equation} \label{eq:genrs}
\begin{split}
Z(t)= \Re \left(2e^{-i\theta(t)}\sum_{n=1}^{n_1} \frac{e^{it \log n}}{\sqrt{n}}\right) + \frac{(-1)^{n_1+1}}{a^{1/2}} &\Phi(\rho)+\\
&\frac{(-1)^{n_1+2}}{96\pi^2 a^{3/2}}\Phi^{(3)}(\rho)+R(t)\,,
\end{split}
\end{equation}

\noindent
where 

\begin{equation}
\Phi(x):=\frac{\cos 2\pi (x^2-x-1/16)}{\cos 2\pi x}\,,  
\end{equation}

\noindent
and $\Phi^{(3)}(x)$ is the third derivative of $\Phi(x)$ with respect to $x$. Gabcke \cite{Ga} showed

\begin{equation}\label{eq:gab}
   |R(t)| \leq .053 t^{-5/4}, \quad \textrm{for $t\geq 200$}\,,
\end{equation}

\noindent
which is sufficient for most applications. Odlyzko and Sch\"onhage \cite{OS} showed how to compute the rotation factor $e^{-i\theta(t)}$, and the correction terms $\Phi(x)$ and $\Phi^{(3)}(x)$, to within $\pm t^{-\kappa-10}$, with $\kappa>0$ fixed, using $t^{o_{\kappa}(1)}$ operations on numbers of $O_{\kappa}(\log t)$ bits. Note throughout, asymptotic constants are taken as $t\to\infty\,$, and the notations $O_{\kappa}(.)$ and $o_{\kappa}(.)$ mean asymptotic constants depend on the parameter $\kappa$ only. Therefore, to calculate the rotated zeta function $Z(t)$, the bulk of the computational effort is in computing the sum 

\begin{equation}\label{eq:msss}
F(t):=\sum_{n=1}^{n_1} \frac{e^{it\log n}}{\sqrt{n}}\,,\qquad n_1=\lfloor \sqrt{t/(2\pi)} \rfloor\,.
\end{equation}

By taking more correction terms in formula (\ref{eq:genrs}), one can arrange for the remainder term $R(t)$ to be bounded by $O(t^{-\kappa})$ for any fixed $\kappa>0$. Odlyzko and Sch\"onhage \cite{OS} showed the additional correction terms can also be computed to within $\pm t^{-\kappa}$ in $t^{o_{\kappa}(1)}$ operations on numbers of $O_{\kappa}(\log t)$ bits. 

Frequently, one is interested in numerically evaluating $Z(t)$ at $N$ points in an interval of the form $t\in [T,T+T^{\alpha}]$, where $\alpha\in[0,1/2]$ say, and $N$ is large. This is the case, for example, when one attempts to locate real zeros of $Z(t)$ to within $\pm T^{-\kappa}$, or study moments of the zeta function. A straightforward application of the Riemann-Siegel formula can do so  in $N\,T^{1/2+o_{\kappa}(1)}$ operations.

The purpose of this note is to improve the running time in the $N$-aspect. So, although the proposed method still consumes $t^{1/2+o_{\kappa}(1)}$ time to evaluate \mbox{$Z(t)$} to within $\pm t^{-\kappa}$ at a {\em{single}} point, it achieves substantially lower running times if one is interested in computing zeta at {\em{many}} points. This type of idea is not new: in the context of the Riemann zeta function, it dates back the algorithm of Odlyzko and Sch\"onhage \cite{OS}.

The main feature of the proposed method is it completely avoids the two essential components of the Odlyzko-Sch\"onhage algorithm, which are the fast Fourier transform, and a sophisticated rational function evaluation algorithm. Instead, the method relies on a straightforward subdivision of the main sum in the Riemann-Siegel formula, a band-limited interpolation technique (see Appendix), and a direct evaluation to obtain the precomputation data. Therefore, its implementation is relatively straightforward, with friendly asymptotic constants, and its error terms are easy to bound. Also, it does not require large amounts of storage space for its precomputation data.  The method far from competes with the Odlyzko-Sch\"onhage algorithm in general. But in many situations, it achieves a similar complexity. 

The basic idea of the method is computing $Z(t)$ for a lot of different, but neighboring, values of $t$ involves many common steps. The method takes advantage of this to achieve lower running times.

Before discussing the method any further, we make a few remarks. By {\em{computing}} zeta we mean to numerically evaluate $Z(t)$ with ``polynomial accuracy,'' that is, with an absolute error bounded by $t^{-\kappa}$, for any fixed $\kappa>0$. We measure the computational complexity of our method by the number of arithmetic operations required: additions, multiplications, divisions, complex exponential, and logarithm (involving numbers of $O_{\kappa}(\log t)$ bits). That in turn can be routinely bounded by the number of bit operations. Lastly, as in the Odlyzko-Sch\"onhage algorithm, the proposed method will generalize easily to any Dirichlet series:

\begin{equation}
\sum_n \frac{a_n}{n^s}\,,\qquad s=\sigma+it\,  \textrm{ with }\, \sigma\in [0,1]\,\textrm{ say}\,, 
\end{equation}

\noindent
assuming the coefficients $a_n$ are known, or can be computed quickly. For definiteness, in the remainder of the paper, we specialize to the rotated zeta function on the critical line $Z(t)$. Our main result is the following, 

\begin{theorem} \label{mainlem}
Given any fixed numbers $\alpha\in[0,1/2]$ and $\kappa>0$, there exists an algorithm that for every $T>100$ will perform $T^{1/2+o_{\kappa}(1)}$ operations on numbers of $O_{\kappa}(\log T)$ bits using $T^{\alpha+o_{\kappa}(1)}$ bits of storage, after which the algorithm will be capable of computing $Z(t)$ at any $t\in [T,T+T^{\alpha}]$ to within $\pm T^{-\kappa}$ in $T^{\alpha+o_{\kappa}(1)}$ operations. \end{theorem}

It is useful to compare the the algorithm of Theorem~\ref{mainlem} with the Odlyzko-Sch\"onhage algorithm (the statement below is equivalent to Theorem 5.1 in \cite{OS} specialized to $Z(t)$):

\begin{os}
Given any $a\in[0,1/2]$, and any constants $\epsilon$ and $\kappa$, there is an effectively computable constant $B=B(\epsilon,\kappa,a)$, and an algorithm that for every $T>0$ will perform $\le B T^{1/2+\epsilon}$ operations on numbers of $\le B\log T$ bits using $\le B T^{a+\epsilon}$ bits of storage and will then be capable of computing any value of $Z(t)$ for $T\le t\le T+T^a$ to within $\pm T^{-\kappa}$ in $\le B T^{\epsilon}$ operations using the precomputed values. 
\end{os} 

As mentioned earlier, the two central ingredients of the Odlyzko-Sch\"onhage algorithm are the fast Fourier transform, and a rational function evaluation algorithm.  Our method completely avoids these central  ingredients. Instead, it relies on a straightforward subdivision of the main sum in the Riemann-Siegel formula, a band-limited interpolation technique (see Appendix), and a direct evaluation for the precomputation data. So, it is significantly simpler.

We remark the Odlyzko-Sch\"onhage algorithm has been implemented at least twice, by Odlyzko \cite{O}, and by Gourdon \cite{G}. Gourdon's implementation replaces the Odlyzko-Sch\"onhage rational function algorithm by the Greengard-Rokhlin algorithm. This was suggested by Odlyzko as a possible improvement in \cite{O}. 

We carry out a comparison between the Odlyzko-Sch\"onhage algorithm and the algorithm of Theorem~\ref{mainlem} in the following context. Suppose we wish to evaluate $Z(t)$ to within $\pm T^{-\kappa}$ at about $\lfloor T^{\alpha}\rfloor$ points in the interval $[T,T+T^{\alpha}]$. This is often the case in applications; for example, when one attempts to locate, to within $\pm T^{-\kappa}$, the ordinates of non-trivial zeros of zeta on the critical line (notice these are the real zeros of $Z(t)$), one expects to require about $T^{o_{\kappa}(1)}$ evaluations of $Z(t)$ per zero. Since there are $T^{\alpha+o(1)}$ zeros of $Z(t)$ in the interval $[T,T+T^{\alpha}]$, then about $T^{\alpha+o_{\kappa}(1)}$ evaluations of $Z(t)$ are needed in total. So, consider Table~\ref{compos}. It compares the running times of the algorithm of Theorem~\ref{mainlem} and the Odlyzko-Sch\"onhage algorithm in such a situation.

\begin{table}[ht]
\renewcommand\arraystretch{1.5}
\footnotesize
\caption{}\label{compos}
\begin{tabular}{l|llll}
Algorithm & Precomputation & Storage & Single Eval. & Average  \\
\hline
Theorem~\ref{mainlem} & $T^{1/2}$ & $T^{\alpha}$ & $T^{\alpha}$ & $T^{1/2-\alpha} + T^{\alpha}$ \\
Odlyzko-Sch\"onhage & $T^{1/2+\epsilon}$ & $T^{\alpha+\epsilon}$ & $T^{\epsilon}$ & $T^{1/2-\alpha+\epsilon}+T^{\epsilon}$ \\
\end{tabular} 
\end{table}

In Table~\ref{compos}, the column ``Single Eval.'' refers to the cost of a single evaluation of $Z(t)$ after the precomputation has been carried out. The column ``Average'' refers to the average cost of evaluating $Z(t)$ at all $\lfloor T^{\alpha}\rfloor$ points (in particular, the column ``Average'' takes the precomputation cost into account). To avoid notational clutter, the complexities listed in the table ignore little-$o$ terms of the form $T^{o(1)}$ (In \cite{OS}, it is stated the $T^{\epsilon}$ term can be replaced by a fixed power of $\log T$; nevertheless, it is included in the table to be consistent with the formal statement of the Odlyzko-Sch\"onhage algorithm). 

In view of the last column of the table, we see if $\alpha\le 1/4$, the algorithm of Theorem~\ref{mainlem} achieves the same amortized complexity as the Odlyzko-Sch\"onhage algorithm, but when $\alpha>1/4$ it does not. This suggests one should restrict $\alpha\in[0,1/4]$. More generally, we have the following corollary to Theorem~\ref{mainlem},

\begin{corollary}\label{corz}
Given any fixed numbers $a\in[0,1/2]$ and $\kappa>0$, there exists an algorithm that for every $T>100$, will be capable of computing $Z(t)$ to within $\pm T^{-\kappa}$ at any $\lfloor T^a \rfloor$ points in $[T,T+T^a]$ using an average of 
  
\begin{displaymath}
\left\{\begin{array}{ll}
T^{1/2-a+o_{\kappa}(1)} & \textrm{operations per point if $0\le a \le 1/4$,}\\
T^{1/4+o_{\kappa}(1)} & \textrm{operations per point if $1/4<a \le 1/2$,}
\end{array}\right.
\end{displaymath}

\noindent
where the operations are performed on numbers of $O_{\kappa}(\log T)$ bits. The storage space requirement for the algorithm is

\begin{displaymath}
\left\{\begin{array}{ll}
T^{a+o_{\kappa}(1)} & \textrm{bits if $0\le a \le 1/4$,}\\
T^{1/4+o_{\kappa}(1)} & \textrm{bits if $1/4<a \le 1/2$.}
\end{array}\right.
\end{displaymath}
\end{corollary}

By comparison, to accomplish the same task as in the statement of the corollary, the Odlyzko-Sch\"onhage algorithm requires $T^{1/2+\epsilon-a+o_{\epsilon,\kappa,a}(1)}$ operations per point on average, and it requires a storage space of $T^{a+\epsilon+o_{\epsilon,\kappa,a}(1)}$  bits (which is the same as corollary~\ref{corz} when $a\in[0,1/4]$). 

The statement of corollary~\ref{corz} follows from a straightforward optimization procedure. Specifically, for any fixed numbers $\delta\ge 0$, $\alpha\in[0,1/2]$, and $\kappa>0$, such that $\delta+\alpha\in[0,1/2]$, consider the successive intervals 

\begin{equation}
[T+jT^{\alpha},T+(j+1)T^{\alpha}]\,, \qquad j=0,1,\ldots,\lceil T^{\delta} \rceil\,.  
\end{equation}

\noindent
Then, applying the algorithm of Theorem~\ref{mainlem} a total of $\lceil T^{\delta} \rceil+1$ times to these intervals, we deduce $Z(t)$ can be computed to within $\pm T^{-\kappa}$ at any single point in

\begin{equation}
[T,T+T^{\alpha+\delta}]\,,\qquad \alpha \in [0,1/2]\,,\qquad \delta \ge 0\,,\qquad \alpha+\delta\in[0,1/2]\,, 
\end{equation}

\noindent
using $T^{\alpha+o_{\kappa}(1)}$ operations, on numbers of $O_{\kappa}(\log T)$ bits, provided a precomputation costing $T^{1/2+\delta+o_{\kappa}(1)}$ operations is performed. Notice the storage space requirement for the precomputation data can always be kept at $T^{\alpha+o_{\kappa}(1)}$ bits. This is because the method will deal with one subinterval $[T+jT^{\alpha},T+(j+1)T^{\alpha}]$ at a time, so the precomputation data for that subinterval can be discarded when the method is done there.

Now, as in the statement of the corollary, let $a\in [0,1/2]$, and suppose we wish to compute $Z(t)$ to within $\pm T^{-\kappa}$ at $\lfloor T^a \rfloor$ points in the interval $[T,T+T^a]$. Optimizing $\alpha$ and $\delta$ to the above situation, we solve 

\begin{equation}
1/2+\delta=a+\alpha\,,\qquad \,\textrm{and}\, \qquad\, \alpha+\delta=a\,,  
\end{equation}

\noindent
which has the solution $\alpha=1/4$ and $\delta=a-1/4$. If $a< 1/4$, then $\delta$ becomes negative. So in this case, we choose $\alpha=a$ and $\delta=0$. For example, when $a\in[0,1/4)$, the method requires $T^{1/2+\delta+o_{\kappa}(1)}=T^{1/2+o_{\kappa}(1)}$ operations on numbers of $O_{\kappa}(\log T)$ bits to compute $Z(t)$ to within $\pm T^{-\kappa}$ at all the $\lfloor T^a\rfloor$ points in the interval $[T,T+T^a]$. This amounts to an average of $T^{1/2-a+o_{\kappa}(1)}$ operations per point. The storage space requirements (for the precomputation data) is $T^{\alpha+o_{\kappa}(1)}=T^{a+o_{\kappa}(1)}$ bits. The analysis of the case $a\in(1/4,1/2]$ is analogous.

We compare the efficiency of our proposed method with Gourdon's implementation \cite{G} of the Odlyzko-Sch\"onhage algorithm. To this end, consider that in order to locate zeros near $t=10^{16}$ to within $\pm 10^{-8}$ say, one expects $\approx 8$ evaluations of $Z(t)$ are required per zero (see \cite{O} p.80, and \cite{G} p.26). The mean spacing of zeros near $10^{16}$ is about $2\pi/\log(10^{16}/(2\pi))\approx 0.18$. Therefore, one expects to evaluate $Z(t)$ at points with $1/8\times 0.18 \approx 0.02$ mean spacing. 

In particular, to locate $2\times 10^9$ consecutive zeros near $t=10^{16}$ one expects to compute $Z(t)$ at $8\times 2\times 10^9=1.6\times 10^{10}$ successive points with mean spacing $0.02$. Results of timing tests (see last entry in Table~\ref{t0} and footnote 3 in Section 3), suggest our method will consume about $9$ minutes to compute $Z(t)$ at $10^5$ such points. By extrapolation, we expect the method to consume $(9/60)\times (1.6\times 10^5)=24,000$ hours to compute $Z(t)$ at $1.6\times 10^{10}$ such points. So, we extrapolate that the method will require 24,000 hours to locate $2\times 10^9$ zeros of $Z(t)$ near $t=10^{16}$. By comparison, Gourdon's implementation of the Odlyzko-Sch\"onhage algorithm (see \cite{G}, p.28) consumes $49.5$ hours to locate the same number of zeros at that height. Thus, Gourdon's implementation is approximately $485$ times faster than our method for that task. In turn, our method is approximately $437$ times faster than $lcalc$'s direct Riemann-Siegel formula; see Table 1, Section 3.

However, if one is interested in finding a smaller set of zeros, then our method might be more suitable, both in terms of timings and the required programming effort. This is because one does not expect the time requirement of the Odlyzko-Sch\"onhage algorithm to decrease significantly as the size of the set of zeros to be located decreases, whereas the time requirement of our method becomes substantially less. The reason is near $t=10^{16}$, and with $10^9$ zeros to be located, we are working in the region $a\approx 1/2$ in corollary~\ref{corz}. So the running time of our method is roughly linear in the number of zeros to be found. For example, to locate $2\times 10^8$ zeros to within $\pm 10^{-8}$ near $t=10^{16}$, we expect our method to consume $10$ times less than in the previous scenario, or $\approx 2,400$ hours. 

We remark our implementation of the method will be available in the next release of Michael Rubinstein's $lcalc$; see \cite{R1}. 

\section{Proof of Theorem~\ref{mainlem}}

In view of formula (\ref{eq:genrs}), and the remarks made thereafter, in order to compute $Z(t)$ with polynomial accuracy, it suffices to numerically evaluate the main sum 

\begin{equation} \label{eq:genrs0}
F(t):=\sum_{n=1}^{n_1} \frac{e^{it\log n}}{\sqrt{n}}\,,\qquad \textrm{where } n_1:=\lfloor \sqrt{t/(2\pi)}\rfloor\,, 
\end{equation}

\noindent
with polynomial accuracy. Given $T>100$, and fixed numbers $\alpha\in[0,1/2]$, $\kappa>0$, we wish to evaluate $F(t)$ for many values of $t\in[T,T+T^{\alpha}]$ to within $\pm T^{-\kappa}$. Let $n_2:=\lfloor \sqrt{T/(2\pi)} \rfloor$. Since $\alpha\in[0,1/2]$, then $n_1$, which is the length of the main sum (\ref{eq:genrs0}), is equal to either $n_2$ or $n_2+1$. Without loss of generality, we may assume the length of the main sum is $n_2$. Let $M:=\min\left\{\lceil T^{\alpha} \rceil\,,\,n_2\right\}$. Note $M$ and $n_2$ depend only on $T$. Then, there exist unique integers $M_1,m\ge 0$ (also depending only on $T$), with $n_2=2^m M+M_1$, and $M_1< 2^m M$. So the main sum can be written as

\begin{equation}\label{eq:initsm} 
\begin{split}
\underbrace{\sum_{1\le n< M} \frac{e^{it\log n}}{\sqrt{n}}}_{\textrm{Initial sum}}+&\sum_{M\le n< 2M} \frac{e^{it\log n}}{\sqrt{n}}+\cdots\\
&+\sum_{2^{m-1}M\le n<  2^m M} \frac{e^{it\log n}}{\sqrt{n}}+\underbrace{\sum_{2^m M\le n\le  2^mM+M_1} \frac{e^{it\log n}}{\sqrt{n}}}_{\textrm{Tail sum}}  
\end{split}
\end{equation}

\noindent
The ``initial sum'' in (\ref{eq:initsm}) can be evaluated directly to within $\pm T^{-\kappa-1}$ in $10M\le 10T^{\alpha}$ operations on numbers of $\lceil (10\kappa+10) \log T \rceil$ bits, say. Thus, we may focus our attention on computing (to within $T^{-\kappa-1}$) the subsums: 

\begin{equation} \label{eq:tailbook}
\sum_{n=2^l M}^{2^{l+1}M-1} \frac{e^{it \log n}}{\sqrt{n}}\qquad l\in\{0,1,\ldots,m-1\}\,,\qquad \textrm{and}\qquad \underbrace{\sum_{n=2^mM}^{2^mM+M_1} \frac{e^{it\log n}}{\sqrt{n}}}_{\textrm{Tail sum}}\,.
\end{equation}

A direct evaluation of all of the subsums (\ref{eq:tailbook}) requires $T^{1/2+o_k(1)}$ operations. But computing the individual terms in these subsums involves many common steps for a lot of different choices of $t\in [T,T+T^{\alpha}]$. We take advantage of this to obtain a substantially lower running times.  So consider one of the subsums on the left side of (\ref{eq:tailbook}). Split it into consecutive blocks of length $2^l$. This gives,

\begin{equation} \label{eq:blockbook}
\begin{split}
\sum_{n=2^lM}^{2^{l+1}M-1} \frac{e^{it \log n}}{\sqrt{n}}=&\underbrace{\sum_{n=2^l M}^{2^l M+2^l-1} \frac{e^{it \log n}}{\sqrt{n}}}_{\textrm{First block}}+\underbrace{\sum_{n=2^l M+2^l}^{2^l M+2(2^l)-1}\frac{e^{it \log n}}{\sqrt{n}}}_{\textrm{Second block}}+ \\
&\underbrace{\sum_{n=2^l M+2(2^l)}^{2^l M+3(2^l)-1}\frac{e^{it \log n}}{\sqrt{n}}}_{\textrm{Third block}}+\cdots +\underbrace{\sum_{n=2^lM+(M-1)2^l}^{2^lM+M2^l-1} \frac{e^{it \log n}}{\sqrt{n}}}_{\textrm{$M^{th}$ block}}\,.
\end{split} 
\end{equation}

\noindent
Notice for any $l\in \{0,1,\ldots,m-1\}$, we have a total of $2^lM/2^l= M$ blocks. So the total number of blocks for any such $l$ is exactly $M$. As for the ``Tail sum'' in (\ref{eq:tailbook}), we use blocks of length $2^m$. This gives a similar outcome as in (\ref{eq:blockbook}), except now the total number of blocks is $M'=\lfloor (M_1+1) 2^{-m}\rfloor$, which is at most $M$, and there is possibly a ``Remainder block'' of length $\le 2^m$,
 
\begin{equation} 
\mathcal{R}_{M,M_1,m}(t):=\underbrace{\sum_{n=2^m M+M'2^m}^{2^mM+M_1} \frac{e^{it\log n}}{\sqrt{n}}}_{\textrm{Remainder block}}\,,\qquad \textrm{where } M':=\left\lfloor \frac{M_1+1}{2^m} \right\rfloor\le M\,.  
\end{equation}

\noindent
Other than the remainder block $\mathcal{R}_{M,M_1,m}(t)$, all our blocks are members of the set 

\begin{equation}
\sum_{n=2^lM+r2^l}^{2^lM+(r+1)2^l-1} \frac{e^{it\log n}}{\sqrt{n}}\,,\qquad l\in \{0,1,\ldots,m\}\,,\,\,\, r\in\{0,1,\ldots,M-1\}\,. 
\end{equation}

\noindent
(when $l=m$, only the values $r< M'$ are relevant to the algorithm). Let 

\begin{equation}
v:=v_{l,r}=2^lM+r2^l\,,\qquad K:=K_l=2^l\,. 
\end{equation}

\noindent
Then the $v^{th}$ block (the one starting at $v:=v_{l,r}$) can be written in the form

\begin{equation} \label{eq:taybook}
\sum_{k=0}^{K-1} \frac{e^{it \log (v+k)}}{\sqrt{v+k}}=e^{it \log v} \sum_{k=0}^{K-1} \frac{e^{it \log (1+k/v)}}{\sqrt{v+k}}=:e^{it\log v} F_{v,K}(t)\,.
\end{equation}

The function $F_{v,K}(t)$ is a {\em{band-limited function}}; that is, its spectrum is limited to a finite interval, which in this case is the interval \mbox{$[0,\log(1+K/v)]$}. In more technical terms, a  band-limited function can be defined as a function whose Fourier transform is a tempered distribution with compact support. There exists a clever method to compute band-limited functions, which we reproduce in the Appendix with slight modifications (see \cite{O}, p.88, for an in-depth discussion and history of the method). To apply band-limited interpolation, first note by construction,

\begin{equation} \label{eq:bdd1}
\frac{K}{v}=\frac{K_l}{v_{l,r}}=\frac{2^l}{2^lM+r2^l}\le \frac{1}{M}\,,
\end{equation} 

\noindent
for all $l\in\{0,1,\ldots,m\}$, and $r\in\{0,1,\ldots,M-1\}$. Therefore, the spectrum of the functions $F_{v,K}(t)$ is always contained in the interval $[0,M^{-1}]$. As for the remainder block $\mathcal{R}_{M,M_1,m}(t)$, let 
\begin{equation}
v':=2^mM+M'2^m\,, \qquad K':=M_1-M'2^m\,, \qquad \textrm{(note $K'\le 2^m$),}  
\end{equation}

\noindent
then write

\begin{equation}
\begin{split} 
\mathcal{R}_{M,M_1,m}(t)&:=\sum_{n=v'}^{v'+K'} \frac{e^{it\log n}}{\sqrt{n}} =e^{it\log v' }\,\sum_{k=0}^{K'} \frac{e^{it\log(1+k/v')}}{\sqrt{v'+k}}\\
&=:e^{it\log v'}\,F_{v',K'}(t)\,.  
\end{split}
\end{equation}

\noindent
Then similarly to the bound (\ref{eq:bdd1}), we have

\begin{equation}
\frac{K'}{v'}=\frac{M_1-M'2^m}{2^mM+M'2^m}\le \frac{2^m}{2^mM+M'2^m}\le \frac{1}{M}\,. 
\end{equation}

\noindent
In particular, the spectrum of $F_{v',K'}(t)$ is also contained in $[0,M^{-1}]$. So, we may apply formulas (\ref{eq:blfi}) and (\ref{eq:kernel}) in the Appendix with $G(t)=F_{v,K}(t)$, and $G(t)=F_{v',K'}(t)$, and (in both cases), 

\begin{equation}
\tau = M^{-1}\,,\,\,\,\,\,\, \beta=3\tau\,,\,\,\,\,\,\,\lambda=(\beta+\tau)/2=2\tau\,,\,\,\,\,\,\, \epsilon_1 =(\beta-\tau)/2=\tau\,,  
\end{equation}

\noindent
say. We then appeal to the bounds (\ref{eq:blfi_error1}) and (\ref{eq:blfi_error}) in the Appendix to obtain

\begin{equation}\label{eq:ttt11}
\begin{split}
F_{v,K}(t) = \frac{\lambda}{\beta} \sum_{|n\pi/\beta-t|<c/\epsilon_1} F_{v,K} \left( \frac{n \pi}{\beta} \right)  \frac{\sin \lambda (n \pi / \beta -t)}{ \lambda ( n \pi / \beta -t)} h ( n \pi &/ \beta -t )\\
&+\mathcal{E}_{v,K}\,.   
\end{split}
\end{equation}

\noindent
for any $c>1$, where

\begin{equation} 
h(u):= \frac{c}{\sinh (c)}\frac{\sinh \sqrt{c^2 - \epsilon_1^2 u^2}}{\sqrt{c^2 - \epsilon_1^2 u^2}}\,,\qquad |\mathcal{E}_{v,K}| < 6 e^{-c}\sum_{k=0}^{K-1} \frac{1}{\sqrt{v+k}}\,.   
\end{equation}

\noindent
By similar calculations, we obtain an analogous formula to (\ref{eq:ttt11}) for $F_{v',K'}(t)$, which corresponds to the remainder block $\mathcal{R}_{M,M_1,m}(t)$. Now, define

\begin{equation}\label{eq:blfi_app}
\tilde{F}_{v,K}(t) = \frac{\lambda}{\beta}  \sum_{|n\pi/\beta-t|<c/\epsilon_1} F_{v,K} \left( \frac{n \pi}{\beta} \right) \frac{\sin \lambda (n \pi / \beta -t)}{ \lambda ( n \pi / \beta -t)} h ( n \pi / \beta -t )\,.  
\end{equation}

\noindent
Also define $\tilde{F}_{v',K'}(t)$ in an analogous way to (\ref{eq:blfi_app}). Note $\tilde{F}_{v,K}$ (also, $\tilde{F}_{v',K'}(t)$) is a sum of $\le 2c\,\beta/\epsilon_1\le 6c$ terms. Then put together, we have

\begin{equation} \label{eq:fin1}
\begin{split}
\underbrace{\sum_{1\le n\le \sqrt{T/(2\pi)}} \frac{e^{it\log n}}{\sqrt{n}}}_{\textrm{Main sum}}=&\underbrace{\sum_{1\le n< M} \frac{e^{it\log n}}{\sqrt{n}}}_{\textrm{Initial sum}}+\sum_{l=0}^{m-1} \sum_{r=0}^{M-1} e^{it \log v_{l,r}} \tilde{F}_{v_{l,r},K_l}(t) \\
&+\sum_{r=0}^{M'-1} e^{it\log v_{m,r}} \tilde{F}_{v_{m,r},K_m}(t) +\tilde{F}_{v',K'}(t)+\mathcal{E}\,,
\end{split}
\end{equation}

\noindent
where,

\begin{equation} \label{eq:fin3}
|\mathcal{E}|\le \sum_{l=0}^{m-1} \sum_{r=0}^{M-1} |\mathcal{E}_{v_{l,r},K_l}|+\sum_{r=0}^{M'-1} |\mathcal{E}_{v_{l,r},K_l}|+ |\mathcal{E}_{v',K'}|\le 20 e^{-c} \sqrt{T}\, 
\end{equation}

\noindent
We choose $c=(\kappa+2)\log T$, so that $|\mathcal{E}|<T^{-\kappa-1}$. Lastly, as we will soon explain, the right side of (\ref{eq:fin1}) can now be evaluated to within $\pm T^{-\kappa}$ using $T^{\alpha+o_{\kappa}(1)}$ operations provided we precompute the entries of the following tables (to within $\pm (6c)^{-1} T^{-\kappa-1}$ each):

\begin{align}\label{eq:tablescount}
&\left\{ F_{v_{l,r},K_l}\left(\frac{n\pi}{\beta}\right)\,:\, l\in[0,m-1],\,r\in[0,M-1]\,,\,\frac{n\pi}{\beta} \in [T-2c/\epsilon_1\,,\,T+T^{\alpha}+2c/\epsilon_1]\right\}\,.\nonumber\\
&\left\{ F_{v_{m,r},K_m}\left(\frac{n\pi}{\beta}\right)\,:\, r\in[0,M'-1]\,,\,\frac{n\pi}{\beta} \in [T-2c/\epsilon_1\,,\,T+T^{\alpha}+2c/\epsilon_1]\right\}\,. \nonumber\\
&\left\{ F_{v',K'}\left(\frac{n\pi}{\beta}\right)\,:\, \frac{n\pi}{\beta} \in [T-2c/\epsilon_1\,,\,T+T^{\alpha}+2c/\epsilon_1]\right\}\,.\nonumber\\ 
\end{align}

\noindent
To compute the entries in the first table to within $\pm (6c)^{-1} T^{-\kappa-1}$ requires a number of operations of at most

\begin{equation}\label{eq:countops}
\begin{split}
\sum_{l=0}^{m-1}\sum_{r=0}^{M-1} &\left\lceil 10 \beta (T^{\alpha}+4c/\epsilon_1)\right\rceil K_l\\
&\le \sum_{l=0}^{m-1}\sum_{r=0}^{M-1} 10 (10+12c)2^l\le 150c\,2^m\,M\le 200c\,T^{1/2}\,, 
\end{split}
\end{equation}

\noindent
where the operations are carried out on numbers of $\lceil (10\kappa+10)\log T\rceil$ bits, say. Note the second inequality in (\ref{eq:countops}) used the fact $\beta/\epsilon_1=3$, which is true by construction. Similarly, the second and third tables in (\ref{eq:tablescount}) require at most $200c\,2^m\,M'$ and $200c\,2^m$ such operations, respectively. Since $M'\le M$, and $2^m M\le T^{1/2}$, then the total cost of precomputing the entries of all three tables in (\ref{eq:tablescount}) is less than

\begin{equation}
600c\,T^{1/2}\le 600(\kappa+2)\, T^{1/2}\log T
\end{equation}

\noindent
operations. Once the precomputation is done, the right side of (\ref{eq:fin1}) can be computed, with the aid of formula (\ref{eq:blfi_app}), to within $\pm T^{-\kappa-1}$ in  less than

\begin{equation}
\begin{split}
 20M+20mM \lceil 2\beta c/\epsilon_1\rceil+20M' \lceil 2\beta c/\epsilon_1 \rceil&+20 \lceil 2\beta c/\epsilon_1\rceil\\
&\le 1000(\kappa+2)\, T^{\alpha} (\log T)^2 
\end{split}
\end{equation}  

\noindent
operations. Finally, the storage space requirement for precomputing the three tables in (\ref{eq:tablescount}) is at most

\begin{equation}
3 m M \lceil 2\beta c/\epsilon_1\rceil \lceil (\kappa+1)\log T\rceil \le 1000(\kappa+2)^2\, T^{\alpha} (\log T)^3 
\end{equation}

\noindent
bits.

\section{Comparison with the Riemann-Siegel formula}

In this section, our algorithm is denoted by BLFI (band-limited function interpolation), and the Riemann-Siegel formula is denoted by RS.

The performance of the BLFI algorithm is compared with that of a relatively optimized version of the RS formula included in the $lcalc$ library. The $lcalc$ library is a C$++$ software developed by Michael Rubinstein to compute values of $L$-functions, including the zeta function.  It can be downloaded at \cite{R1}.  

The BLFI algorithm was also coded in C$++$, essentially as described in the previous sections\footnote{Our implementation of the BLFI algorithm differed from the description in the previous sections only in that we centered the band-limited functions $F_{v,K}$ so their spectrum is in $[-0.5\,\log(1+K/v),0.5\,\log(1+K/v)]$ rather than $[0,\log(1+K/v)]$. This allows the recovery of values of the functions $F_{v,K}$ using less frequent sampling; see \cite{O}, pp. 90-91.}. We used double-precision arithmetic, which allows a maximum of $16$ digits of accuracy\footnote{Working with $30$-digit arithmetic takes about $10$ times longer for both RS and BLFI. But the slowdown can be made much less severe by following some of the tricks in \cite{G} p.13}.

We used the BLFI algorithm and the RS formula to evaluate $\zeta(1/2+it)$ on a grid of points of the form

\begin{equation}
[T\,,T+n\Delta]\,,\qquad n=1,\ldots,N\,.  
\end{equation}

\noindent
Here, $\Delta>0$ is the spacing (or density) of points\footnote{The grid of points where the zeta function is to be evaluated need not at all be uniform. What matters is the average spacing  of the grid points.}, and $N$ is the number of points in the grid. In case of the BLFI algorithm, an upper bound for the truncation error $\mathcal{E}$ was also chosen. The output of the precomputation needed by the BLFI algorithm is stored in dynamic memory, not saved in files. 

The running times of the BLFI algorithm, presented in tables below, account for everything, including the precomputation. Also, when $T>10^{10}$, timings for the RS formula were obtained by evaluating zeta at a number of points $M<N$, then multiplying the running time by $N/M$. 

In the tables to follow, the running times are formatted as  $x\,\textrm{m } \,y\,\textrm{s }$, where $x$ and $y$ are the numbers of minutes and seconds, respectively, consumed by the method under test. Finally, our tests were carried out on a personal Mac machine. The same compiler options were used for both programs. 

Tables~\ref{t0} and~\ref{t1} list running times for BLFI and RS at various heights. The last column is the ratio of the time consumed by RS to the time consumed by the BLFI algorithm. At each height, we used a grid of $10^5$ equidistant points. In Table~\ref{t0}, the spacing of the grid points is $\Delta=0.01$, and in Table~\ref{t1} it is $\Delta=0.1$. The  truncation error was chosen to satisfy $\mathcal{E}<10^{-8}$ in both tables. 

\begin{table}[ht]
\footnotesize
\caption{\footnotesize Running times of BLFI and RS with $N=10^5$, $\Delta=0.01$, and $\mathcal{E}<10^{-8}$.}\label{t0}
\renewcommand\arraystretch{1.5}
\begin{tabular}{l|l|l|l}
$T$       & RS        & BLFI      & Ratio \\
\hline
$10^8$    & 0m 18s  &  0m 7s & 2  \\
$10^{10}$ & 2m 56s  &  0m 12s & 14   \\
$10^{12}$ & 29m 10s  & 0m 35s  & 50   \\
$10^{14}$ & 348m 0s  & 2m 35s  & 134   \\
$10^{16}$ & 3700m 0s  & 8m 28s  & 437 \\  
\end{tabular}
\end{table}

\begin{table}[ht]
\footnotesize
\caption{\footnotesize Running times of BLFI and RS with $N=10^5$, $\Delta=0.1$, and $\mathcal{E}<10^{-8}$.}\label{t1}
\renewcommand\arraystretch{1.5}
\begin{tabular}{l|l|l|l}
$T$       & RS        & BLFI      & Ratio \\
\hline
$10^8$    & 0m 18s  & 0m 9s  & 2   \\
$10^{10}$ & 2m 56s  & 0m 27s  & 6   \\
$10^{12}$ & 29m 10s  & 1m 34s  & 18   \\
$10^{14}$ & 348m 0s  & 5m 50s  & 60   \\
$10^{16}$ & 3700m 0s  & 18m 45s  & 205  
\end{tabular}
\end{table}

Tables~\ref{t0} and~\ref{t1} indicate the running time of RS grows like $T^{1/2}$, as expected, while the running time of the BLFI algorithm grows roughly like $T^{1/4}$, also as expected. As $T$ increases, the savings achieved by the BLFI algorithm are accentuated. For example, when $T=10^{16}$, and $\Delta=0.1$, the BLFI algorithm is more than 200 times faster than RS, but it is only $2$ times faster when $T=10^8$.   

We measure the sensitivity of the running time of the BLFI algorithm to perturbations in the values of its parameters and input. Table~\ref{t2} indicates the running time of the BLFI algorithm grows roughly linearly with the number of grid points $N$, as soon as $N$ is large enough, which is expected. Table~\ref{t3} shows the running time is not radically affected by changes in the error allowance $\mathcal{E}$. This is not surprising either, because demanding higher precision increases the number of terms in the BLFI formula (\ref{eq:ttt11}) only logarithmically. 

\begin{table}[ht]
\footnotesize
\caption{\footnotesize Running times of BLFI with $T=10^{12}$, $\Delta=0.1$, and $\mathcal{E}<10^{-8}$.}\label{t2}
\renewcommand\arraystretch{1.5}
\begin{tabular}{l|l}
$N$      & BLFI     \\
\hline
$10^3$   &  0m 2s\\
$10^4$   &  0m 10s\\
$10^5$   &  1m 34s\\
$10^6$   &  15m 27s
\end{tabular}
\end{table}

\begin{table}[ht]
\footnotesize
\caption{\footnotesize Running times of BLFI with $T=10^{12}$, $\Delta=0.1$, and $N=10^5$.}\label{t3}
\renewcommand\arraystretch{1.5}
\begin{tabular}{l|l}
$\mathcal{E}$ & BLFI \\
\hline
$10^{-6}$  &1m 28s  \\ 
$10^{-8}$  &1m 34s  \\
$10^{-10}$ &1m 41s  \\
$10^{-12}$ &1m 48s  \\
$10^{-14}$ &1m 53s           
\end{tabular}
\end{table}

Lastly, Table~\ref{t4} lists timings for the BLFI Algorithm for various values of $\Delta$. As $\Delta$ increases, the grid expands, so more effort is exerted during the precomputation. This explains the slowdown in the method as $\Delta$ increases.  

\begin{table}[ht]
\footnotesize
\caption{\footnotesize Running times of BLFI with $T=10^{12}$, $N=10^5$, and $\mathcal{E}<10^{-8}$.}\label{t4}
\renewcommand\arraystretch{1.5}
\begin{tabular}{l|l}
$\Delta$  & BLFI     \\
\hline
0.01      & 0m 35s \\
0.05      & 1m 10s \\ 
0.1       & 1m 34s \\
0.2       & 2m 8s \\
0.4       & 2m 54s 
\end{tabular}
\end{table}

\section{Appendix: band-limited interpolation}

Much of the material in this section is contained in \cite{O}, pages 88-93. For the convenience of the reader, it is reproduced here with slight modifications. Let

\begin{equation}\label{eq:eqgt}
G(t) = \int_{- \tau}^\tau  g(x) e^{ixt} dx\,,
\end{equation}

\noindent
be a band-limited function, where $g(x)$ is a (finite) linear combination of delta functions supported on $(-\tau,\tau)$. It is well-known $G$ can be recovered completely from its values at the grid points $\{n\pi/\tau\,|\, n\in \mathbf{Z}\}$. But this recovery process is not efficient, because  it requires the values of $G$ at many sample points $n\pi/\tau$. 

The idea of a band-limited interpolation technique is to sample $G$ on a denser grid, say $\{n\pi/\beta\,|\,n\in\mathbf{Z}\}$, where $\beta>\tau$, after which $G(t)$ can be recovered much more quickly, from its values at only a few grid points $n\pi/\beta$ that are close to $t$.   

Specifically, choose $\beta >\tau$, and define $\lambda :=( \beta + \tau ) /2$, $\epsilon_1:=(\beta-\tau)/2$. Let $I$  denote the characteristic function of the interval $[-\lambda,\lambda]$, and let $H$ be any continuous function with total mass 1 supported on $[-\epsilon_1, \epsilon_1]$. Define $f:=I*H$, the convolution of $I$ and $H$, and let $\hat{f}$ and $\widehat{H}$ denote the Fourier transforms:

\begin{equation}
\hat{f}(t):=\int_{-\infty}^{\infty} f(x)e^{-itx}\,dx\,,\qquad \widehat{H}(t):= \int_{-\infty}^{\infty} H(x) e^{-itx}\,dx\,. 
\end{equation}

By construction, $f=I*H$ is supported on $[-\beta,\beta]$, and it is identically 1 on $[-\tau,\tau]$. And by hypothesis, $g(x)$ is supported on $[-\tau,\tau]$.  Therefore, $f(x)g(x)\equiv g(x)$, the only difference is the left side involve smoothing in the ``redundant'' interval $[\tau,\beta]\cup [-\beta,-\tau]$. The latter observation is what gives the band-limited technique its edge. This is because some Fourier analysis yields,   

\begin{equation}\label{eq:blfi}
\begin{split}
G(t) &= \int_{-\infty}^{\infty} f(x) g(x) e^{ixt} =\frac{\lambda}{\beta}  \sum_n  G( n \pi / \beta )\,  \hat{f}(t-n\pi/\beta) \\ 
&=  \frac{\lambda}{\beta}  \sum_{n=-\infty}^{\infty}  G( n \pi / \beta )  \frac{\sin \lambda (t- n \pi / \beta)} {\lambda (t-n \pi / \beta)} \widehat{H}(t- n \pi / \beta )\,, 
\end{split}
\end{equation}

\noindent
In particular, we can try to choose the smoothing function $\widehat{H}$ so as to accelerate the convergence of the right side in (\ref{eq:blfi}). There are many choices for $\widehat{H}$ (e.g. a Gaussian). Following Odlyzko \cite{O}, we choose the following function, which still closely resembles a Gaussian, but leads to smaller big-$O$ constants (see \cite{L} and \cite{O}, p.92):  

\begin{equation} \label{eq:kernel}
\widehat{H}(u) = \frac{c}{\sinh (c)}\frac{\sinh \sqrt{c^2 - \epsilon_1^2 u^2}}{\sqrt{c^2 - \epsilon_1^2 u^2}}\,,
\end{equation}

\noindent
where $c>1$ is any fixed number.  Finally, we only sum the terms with $| n \pi / \beta - t | < c/ \epsilon_1$ in formula (\ref{eq:blfi}), which yields a truncation error $\mathcal{E}_1$ satisfying

\begin{equation} \label{eq:blfi_error1} 
|\mathcal{E}_1|< 2 \,\sup_{t\in \mathbb{R}} |G(t)|\int_{|u| > c/\epsilon_1}  \left|\frac{\widehat{H}(u)}{u}\right| du\,, 
\end{equation}

\noindent
This is bounded by (see \cite{O}, p.92, and \cite{L}):

\begin{equation} \label{eq:blfi_error}
|\mathcal{E}_1|< 2\,\sup_{t\in \mathbb{R}} |G(t)|  \log  \frac{1+ e^{-c}}{1-e^{-c}}\le 6\,\sup_{t\in\mathbb{R}}|G(t)| e^{-c}\,.
\end{equation}
\\

\noindent
{\bf{Acknowledgement.}} The author would like to thank Michael Rubinstein for helpful discussions on the topic of this paper. The author would like to acknowledge the many helpful comments by the anonymous referee.  

\bibliographystyle{mcom-l}

\end{document}